%% file: article.tex
\documentclass{gtart}

\newcommand{\pathtotrunk}{./}
\input{article/article_preamble.tex}
\input{article/top_matter.tex}

\begin{document}

\begin{abstract}
Let $V$ be the $7$-dimensional irreducible representation of the quantum
group $U_q(\mathfrak{g}_2)$. For each $n$, there is a map from the braid
group $\cB_n$ to the endomorphism algebra of the $n$-th tensor power of
$V$, given by $\cR$ matrices. Extending linearly to the braid group
algebra, we get a map $$\cA \cB_n \to
\End{U_q(\mathfrak{g}_2)}{V^{\tensor n}}.$$ Lehrer and Zhang
\cite{MR2271576} prove this map is surjective, as a special case of a more general result.

Using Kuperberg's spider for $G_2$ from \cite{MR1403861}, we give an
elementary diagrammatic proof of this result.
\end{abstract}

\maketitle


\section{Kuperberg's spider for $G_2$}
We recall just enough from \cite{MR1265145, MR1403861} for our purposes.

We fix the ground ring $\cA = \Complex(q)$.  Kuperberg's $q$ is $q^2$ here, bringing our conventions into agreement with those for quantum groups as presented in \cite{MR1358358, MR1359532, MR2286123}.

First consider the braided tensor category $\Rep'(G_2)$, with objects
tensor powers of the $7$-dimensional representation of
$U_q(\mathfrak{g}_2)$, and morphisms linear maps commuting with the
actions of $U_q(\mathfrak{g}_2)$. (The prime in the notation indicates
this is just a full subcategory of the actual representation category; we
don't allow arbitrary representations, although every representation does appears as a subobject of some tensor power.)

Second consider the category $\cT(G_2)$, with objects natural numbers,
and morphisms planar trivalent graphs embedded in a rectangle, modulo
certain relations. A morphism from $n$ to $m$ should be a graph with $n$
boundary points along the bottom edge of the rectangle, and $m$ boundary
points along the top edge. Composition of morphisms is vertical stacking.
The category becomes a tensor category by adding natural numbers at the
level of objects, and juxtaposing graphs side by side at the level of
morphisms. The relations from \cite{MR1403861} are shown in Figure
\ref{fig:G2-relations}. The category $\cT(G_2)$ becomes a braided tensor
category with the formulas for a crossing from \cite{MR1403861},
\begin{equation}
\label{eq:crossing}
\mathfig{0.08}{crossing} = \frac{1}{1+q^{-2}} \mathfig{0.052}{I} + \frac{1}{1+q^2} \mathfig{0.096}{H} + \frac{1}{q^{2}+q^{4}} \mathfig{0.08}{cap-cup} + \frac{1}{q^{-2}+q^{-4}} \mathfig{0.08}{id}
\end{equation}

\begin{figure}[!ht]
\begin{align*}
\mathfig{0.04}{loop} &= q^{10} + q^{8} + q^2 + 1 + q^{-2} + q^{-8} + q^{-10} \\
\mathfig{0.04}{lollipop} &= 0 \\
\mathfig{0.04}{bigon} &= -(q^6 + q^4 + q^2 + q^{-2} + q^{-4} + q^{-6}) \mathfig{0.002}{strand}\\
\mathfig{0.048}{triangle} &=  q^4 + 1 +q^{-4} \mathfig{0.048}{vertex}\\
\mathfig{0.044}{square} &=  -(q^2+q^{-2}) \left(\mathfig{0.034}{I} + \mathfig{0.05}{H} \right)+ (q^2+1+q^{-2}) \left(\mathfig{0.034}{cap-cup} + \mathfig{0.034}{id} \right) \\
\mathfig{0.052}{pentagon} &=  \mathfig{0.052}{tree1} + \mathfig{0.052}{tree2} + \mathfig{0.052}{tree3} + \mathfig{0.052}{tree4} + \mathfig{0.052}{tree5} \\
 & \qquad - \mathfig{0.052}{forest1}-\mathfig{0.052}{forest2}-\mathfig{0.052}{forest3}-\mathfig{0.052}{forest4}-\mathfig{0.052}{forest5}\\
\end{align*}
\caption{The relations for the $G_2$ spider, from \cite{MR1265145, MR1403861}
omitting relations involving double edges which we don't use. Note also that there is a sign error in \cite{MR1403861} but not in \cite{MR1265145}.}
\label{fig:G2-relations}
\end{figure}

We now need two theorems from \cite{MR1403861}, one easy, one hard.
\begin{thm}[Web bases for $\cT(G_2)$]
\label{thm:bases}
Each $\operatorname{Hom}$ space in $\cT(G_2)$ has a basis given by
diagrams with no loops, in which each internal face has at least $6$
edges.
\end{thm}
\begin{thm}[Isomorphism]
\label{thm:iso}
The categories $\Rep'(G_2)$ and $\cT(G_2)$ are equivalent as braided
tensor categories.
\end{thm}

\section{Diagrammatic proof}

Our goal is now to prove that an arbitrary trivalent graph in the rectangle (with equal numbers of boundary points along the top and bottom edges) can be written, modulo the relations above, as a $\mathcal{A}$-linear combination of diagrams coming from braids via Equation \eqref{eq:crossing}. Using Kuperberg's Theorems \ref{thm:bases} and \ref{thm:iso}, this would then give a combinatorial proof of Lehrer and Zhang's 
\begin{thm}[Surjectivity]
\label{thm:surjectivity}
The map from the braid group algebra to endomorphisms of tensor powers of the $7$-dimensional irreducible representation of $U_q(\mathfrak{g}_2)$
$$\cA \cB_n \to
\End{U_q(\mathfrak{g}_2)}{V^{\tensor n}}$$
is surjective.
\end{thm}
\begin{rem}
In fact, their result holds for any ``strong multiplicity free'' representation of a quantum group.
\end{rem}

Beginning over-optimistically, we might guess that any such diagram can in fact be written as a composition of factors
\begin{equation*}
\mathfig{0.25}{variable-cap-cup}, \mathfig{0.25}{variable-I} \quad \text{and} \quad \mathfig{0.25}{variable-H}
\end{equation*}
and then make use of the (presumably easy) special case of the Theorem for $n=2$. Even though this is false, it's the right direction and contains the essential idea. Small counterexamples to this guess are provided by
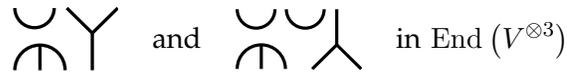
\begin{equation*}
\mathfig{0.1}{down-Y-cup} \quad \text{and} \quad \mathfig{0.12}{down-L-cups} \quad \text{in $\End{}{V^{\tensor 3}}$}
\end{equation*}

The correct argument will involve three steps. First, we'll prove a little lemma allowing us to rearrange connected components of a graph. Second, we'll prove that the image of the braid group does hit a certain finite list of small graphs (including the examples above). Third, we'll use euler measure to inductively rewrite a graph in terms of linear combinations of braids and a graph with fewer vertices, eventually getting down to graphs in our list.

\begin{lem}
\label{lem:tensor}
Suppose a graph $D_1 \tensor D_2$ is a tensor product of two diagrams. Then $D_1 \tensor D_2$ is in the image of the braid group if and only if $D_2 \tensor D_1$ is.
\end{lem}
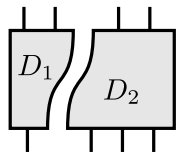
\begin{proof}
If $$\mathfig{0.16}{D1D2}$$ is in the image of the braid group, so is $$\mathfig{0.16}{D1D2braided}.$$ Since the $G_2$ spider is a braided tensor category, this is exactly the same as $$\mathfig{0.16}{D2D1}.$$
\end{proof}

\begin{lem}
\label{lem:braid2}
The graphs
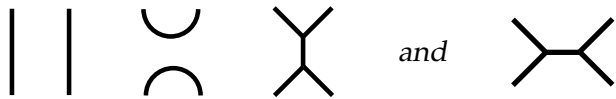
\begin{equation*}
\mathfig{0.06}{id} \qquad \mathfig{0.06}{cap-cup} \qquad \mathfig{0.06}{I} \qquad \text{and} \qquad \mathfig{0.1}{H}
\end{equation*}
are in the image of $\mathcal{AB}_2$. 
\end{lem}
\begin{proof}
These four diagrams are linearly independent and span $\End{}{V^{\tensor 2}}$. Moreover, the braiding lives in the same space, and has four distinct eigenvalues (c.f. \cite{d2n-links}), and so its powers also span $\End{}{V^{\tensor 2}}$. For reference, we'll give the characteristic equation:
\begin{multline*}
\sigma^4 = -q^{-16} \sigma^0 + (q^{-18} - q^{-16} - q^{-10} + q^{-4}) \sigma^1 + \\ + (q^{-18} + q^{-12} - q^{-10} - q^{-6} + q^{-4} + q^2) \sigma^2 + (q^{-12} - q^{-6} - 1 + q^2) \sigma^3.
\end{multline*}
Explicit formulas appear in the appendix, along with instructions for using the short computer program that produces them.
\end{proof}

\begin{lem}
\label{lem:braid3}
The graphs
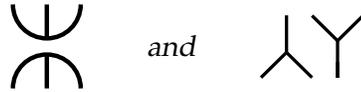
\begin{equation*}
\mathfig{0.07}{down-up} \qquad \text{and} \qquad \mathfig{0.1}{LY}
\end{equation*}
 are in the image of $\mathcal{AB}_3$.
\end{lem}
\begin{proof}
We just do the computations directly, with the help of a computer, obtaining the formulas which appear in the appendix. In particular, we find $35$ braids that provide a basis for $\End{}{V^{\tensor 3}}$.
\end{proof}

\begin{cor}
\label{cor:braid456}
The graphs
\begin{equation*}
\mathfig{0.1}{down-Y-cup}, \mathfig{0.125}{Ls-cup}, \mathfig{0.11}{downs-cups},  \mathfig{0.11}{down-L-cups}  \qquad \text{and} \qquad \mathfig{0.18}{down-Ys}
\end{equation*}
are in the image of the appropriate braid groups.
\end{cor}
\begin{proof}
We'll show that each of these graphs is actually in the braided tensor subcategory generated by the graphs appearing in the two previous lemmas. (In fact, without even needing to take linear combinations.)
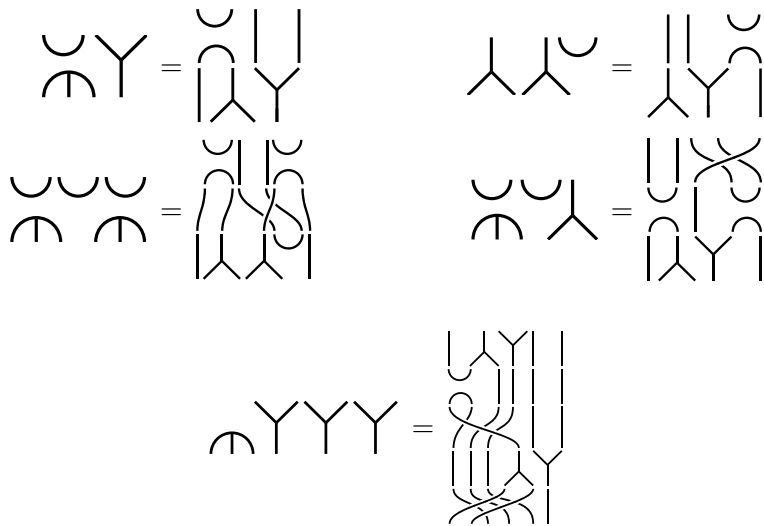
\begin{align*}
\mathfig{0.1}{down-Y-cup} & = \mathfig{0.1}{down-Y-cup-decomposition} &
\mathfig{0.125}{Ls-cup} & =  \mathfig{0.11}{Ls-cup-decomposition} \\
\mathfig{0.13}{downs-cups} & = \mathfig{0.11}{downs-cups-decomposition} &
\mathfig{0.12}{down-L-cups} & = \mathfig{0.11}{down-L-cups-decomposition} 
\end{align*}
\begin{align*}
\mathfig{0.18}{down-Ys} & = \mathfig{0.11}{down-Ys-decomposition} 
\end{align*}
\end{proof}

\begin{prop}[Base case for the induction]
Any $\cT(G_2)$ basis diagram in which every connected component contains at most one vertex is in the image of the braid group.
\end{prop}
\begin{proof}
Any such diagram is just a tensor product of factors, each of which is one of the graphs
\begin{equation*}
\mathfig{0.04}{cup}, \mathfig{0.04}{cap}, \mathfig{0.04}{up}, \mathfig{0.04}{down}, \mathfig{0.04}{Y} \quad \text{or} \quad \mathfig{0.04}{L}.
\end{equation*}
By Lemma \ref{lem:tensor}, we can take this tensor product in any order we like. We claim that any such tensor product is actually a tensor product of the graphs in Lemmas \ref{lem:braid2} and \ref{lem:braid3} and Corollary \ref{cor:braid456}, and so the entire diagram is in the image of the braid group. The argument is a tedious but straightforward case-bash. Suppose we have the diagram
\begin{equation*}
\mathfig{0.03}{cup}^{\tensor a} \mathfig{0.03}{cap}^{\tensor b} \mathfig{0.03}{up}^{\tensor c} \mathfig{0.03}{down}^{\tensor d} \mathfig{0.03}{Y}^{\tensor e} \mathfig{0.03}{L}^{\tensor f} .
\end{equation*}
Note that $2a+3c+e = 2b+3d+f$, by counting boundary points at top and bottom.
By splitting off factors of $\mathfig{0.045}{cap-cup}$, we can assume that at most one of $a$ and $b$ is nonzero. Similarly, splitting off $\mathfig{0.045}{down-up}$ or $\mathfig{0.08}{LY}$, at most one of $c$ and $d$ and at most one of $e$ and $f$ are nonzero. Let's assume without loss of generality that $b=0$.

If $a > 0$, then either we can split off copies of $\mathfig{0.08}{Ls-cup}$, or $f<2$. If $f = 1$, $e$ must be zero, and we have $2a + 3c = 3d + 1$, and so $a \cong 2 \pmod{3}$. Thus $d\geq 1$, and we can split off a copy of $\mathfig{0.08}{down-L-cups}$. Otherwise, if $f=0$, we have $2a + 3c + e = 3d$, so $d \geq 1$ and $c = 0$. Now, since $2a+e=3d \geq 3$, either $a \geq 3$, $e \geq 3$, or $a,e \geq 1$. In each of these cases, we can split something off; either $\mathfig{0.08}{downs-cups}$, $\mathfig{0.1}{down-Ys}$ or $\mathfig{0.07}{down-Y-cup}$ respectively.

On the other hand, if $a=0$, let's further assume without loss of generality that $c=0$. We thus have $e = 3d + f$; since $e = 0$ would imply we have the empty diagram, $f=0$ instead, and the entire diagram is just a tensor power of $\mathfig{0.1}{down-Ys}$.
\end{proof}

\begin{proof}[Proof of Theorem \ref{thm:surjectivity}]
Suppose now we have an arbitrary basis diagram. We will show that a connected component with at least two vertices has either a $\mathfig{0.045}{I}$ or a $\mathfig{0.06}{H}$ attached along its top or bottom edge. Since these two diagrams are Laurent polynomials in the positive crossing, we can rewrite the basis diagram as the product of something in the image of the braid group, and another basis diagram with fewer vertices. Repeating this, we reduce to the case that no connected component contains more than one vertex, at which point we're done by the previous Proposition.

The euler measure argument is straightforward. Consider a connected component of the basis diagram with at least two vertices. Assign formal angles of $\frac{2\pi}{3}$ around trivalent vertices, and of $\frac{\pi}{2}$ on either side of an edge meeting the boundary. (We \emph{don't} assign angles to the corners of the rectangle, so the total euler measure will be the euler measure of the disc, $+1$.) Since internal faces of the component have at least 6 edges, by Theorem \ref{thm:bases}, they have non-positive euler measure. Let $b_k$ be the number of faces adjacent to the boundary, and meeting $k$ edges of the graph. These faces meet the boundary of the disc surrounding the component exactly once, since the component is connected. Thus the euler measure of a face counted by $b_k$ is $$\frac{1}{4} + \frac{1}{4} + \frac{1}{3}(k-1) - \frac{1}{2}(k+1) + 1 = \frac{2}{3} - \frac{k}{6}.$$ The number $b_1$ is zero (this could only occur if the component were a single strand, but it must have at least two vertices), and so we obtain
$$\sum_{k\geq 2} b_k \left(\frac{2}{3} - \frac{k}{6}\right) \geq  1,$$
and so
$$\frac{b_2}{3}+\frac{b_3}{6}\geq 1,$$
which we soften to $b_2 + b_3 \geq 3$.
There are thus at least three faces touching either $2$ or $3$ edges of the component. At least one of these must be attached to the top or bottom of the rectangle, avoiding the sides. If that face touches $3$ edges, we're done, as there must be a $\mathfig{0.06}{H}$ adjacent to the boundary. If it only touches $2$ edges, the hypothesis that the component is connected and has at least two vertices ensures that there's a $\mathfig{0.045}{I}$ adjacent to the boundary.
\end{proof}

\section{Questions}
We end with two questions relating Kuperberg's spider and the category of tilting modules at a root of unity.
\begin{question}
We can specialize Kuperberg's spider for $G_2$ to any root of unity $q$. The braiding is still defined as long as $q+q^{-1} \neq 0$, and the braiding is still surjective as long as 
$$ \left(q^8-1\right) \left(q^4-q^2+1\right)  \left(q^6 + q^4 + q^2 + 1\right) \neq 0.$$
(See the explicit formulas in the appendix.) We can idempotent complete, and quotient by negligibles to produce a semisimple category. Is this category equivalent to the semisimple quotient of the category of tilting modules for the integral form of $U_q(G_2)$?
\end{question}

More optimistically, we might hope for an affirmative answer to:
\begin{question}
Are these categories equivalent even before we quotient by negligible morphisms on either side?
\end{question}

Results in \cite{MR1265145} and \cite{d2n-links} ensure that there is a functor from the specialization of the spider to the category of tilting modules. I'd hoped to be able to leverage the surjectivity of the braiding map into a proof that this functor was an equivalence, but I still don't see how to do this.

\appendix
\section{Explicit formulas}
\newcommand{\code}[1]{{\tt{#1}}}
By downloading the sources for this article from the \code{arXiv}, you'll find a \code{Mathematica} notebook at \code{code/formulas.nb}. This notebook relies on the \code{QuantumGroups`} package, which you can download as part of the \code{KnotTheory`} package from \cite{katlas}. After setting the path in the first line to point at your copy of the \code{QuantumGroups`} package, you can run the remaining lines to produce the formulas appearing below.
\paragraph{Formulas for Lemma \ref{lem:braid2}}
\begin{align*}
\mathfig{0.06}{cap-cup} & = \left((q^2 - 1) (q^4 - q^2 + 1) (q^6 + q^4 + q^2 + 1)\right)^{-1} \times \\ & \qquad \Big(\left(q^{22}\right) \sigma^0 + \left(-q^{20}+q^{22}+q^{28}\right) \sigma^1 + \left(-q^{20}-q^{26}-q^{28}\right)\sigma^2 + \left(-q^{26}\right) \sigma^3 \Big)\\
\mathfig{0.06}{I} & = \left((q^2 - 1) (q^4 - q^2 + 1) (q^6 + q^4 + q^2 + 1)\right)^{-1} \times \\ & \qquad \big(\left(-q^8-q^{12}\right) \sigma^0 + \left(q^6-q^8+q^{10}-q^{12}+q^{20}+q^{24}\right) \sigma^1 + \\ & \qquad + \left(q^6+q^{10}-q^{18}+q^{20}-q^{22}+q^{24}\right)\sigma^2 + \left(-q^{18}-q^{22}\right) \sigma^3 \big)\\
\mathfig{0.08}{H} & = \left((q^2 - 1) (q^4 - q^2 + 1) (q^6 + q^4 + q^2 + 1)\right)^{-1} \times \\ & \qquad \big(\left(q^2-q^{4}+2q^{6}+q^{12}-q^{14}-q^{18}\right) \sigma^0 + \\ & \qquad + \left(-q^{-2}-2q^{4}+2q^{6}-q^{8}+q^{10}+q^{12}+q^{16}-2q^{18}-q^{22}-q^{24}\right) \sigma^1 + \\ & \qquad +  \left(-q^{4}-q^{8}+2q^{16}-q^{18}+q^{20}-q^{24}\right)\sigma^2 + \left(q^{16}+q^{20}+q^{22}\right) \sigma^3 \big)\\
\end{align*}

\paragraph{Formulas for Lemma \ref{lem:braid3}}
\begin{align*}
\mathfig{0.06}{down-up} =&  \left(\left(q^2-1\right)^2 \left(q^2+1\right) \left(q^4-q^2+1\right)^2\right)^{-1} \times \Bigg( \\ &
\scalebox{0.7}{$-q^2 \left(q^2-1\right) \left(q^{26}-q^{24}+q^{22}+q^{20}-2 q^{18}+4 q^{16}-3 q^{14}+2 q^{10}-4 q^8+q^6-q^4-1\right) $}  \mathbf{1}
 \displaybreak[1] \\  & \qquad 
\scalebox{0.7}{$+q^8 \left(-q^{22}+3 q^{20}-5 q^{18}+4 q^{16}-3 q^{12}+7 q^{10}-6 q^8+2 q^6+q^4-4 q^2+2\right) $}  \sigma_{1}
 \displaybreak[1] \\  & \qquad 
\scalebox{0.7}{$-q^4 \left(q^2-1\right) \left(q^{18}-q^{16}+2 q^{14}+2 q^8-q^6+q^4+1\right) $}  \sigma_{1}^{-1}
 \displaybreak[1] \\  & \qquad 
\scalebox{0.7}{$+q^2 \left(q^2-1\right) \left(q^{22}-q^{20}+q^{18}+3 q^{12}-2 q^{10}+3 q^8+q^6-2 q^4+2 q^2-1\right) $}  \sigma_{2}
 \displaybreak[1] \\  & \qquad 
\scalebox{0.7}{$-q^4 \left(q^2-1\right) \left(2 q^{14}+q^{10}+3 q^8-q^6+q^4+1\right) $}  \sigma_{2}^{-1}
 \displaybreak[1] \\  & \qquad 
\scalebox{0.7}{$+q^{10} \left(q^2-1\right)^2 \left(q^2+1\right) \left(q^{12}-q^{10}+3 q^8-q^6+3 q^4+1\right) $}  \sigma_{1}\sigma_{1}
 \displaybreak[1] \\  & \qquad 
\scalebox{0.7}{$+\left(-q^{24}+3 q^{22}-3 q^{20}+q^{18}+q^{16}-5 q^{14}+5 q^{12}-4 q^{10}+q^8+2 q^6-2 q^4+q^2\right) $}  \sigma_{1}\sigma_{2}
 \displaybreak[1] \\  & \qquad 
\scalebox{0.7}{$-q^{10} \left(q^{10}-2 q^8+q^6-q^4-q^2+1\right) $}  \sigma_{1}\sigma_{2}^{-1}
\scalebox{0.7}{$+q^4 \left(q^4-q^2+1\right) \left(q^8+q^2-1\right) $}  \sigma_{1}^{-1}\sigma_{2}
 \displaybreak[1] \\  & \qquad 
\scalebox{0.7}{$-q^6 \left(q^8-q^6+q^4-q^2+1\right) $}  \sigma_{1}^{-1}\sigma_{2}^{-1}
 \displaybreak[1] \\  & \qquad 
\scalebox{0.7}{$+\left(-q^{24}+3 q^{22}-3 q^{20}+q^{18}+q^{16}-5 q^{14}+5 q^{12}-4 q^{10}+q^8+2 q^6-2 q^4+q^2\right) $}  \sigma_{2}\sigma_{1}
 \displaybreak[1] \\  & \qquad 
\scalebox{0.7}{$+q^4 \left(q^4-q^2+1\right) \left(q^8+q^2-1\right) $}  \sigma_{2}\sigma_{1}^{-1}
\scalebox{0.7}{$+\left(q^{22}-q^{24}\right) $}  \sigma_{2}\sigma_{2}
 \displaybreak[1] \\  & \qquad 
\scalebox{0.7}{$-q^{10} \left(q^{10}-2 q^8+q^6-q^4-q^2+1\right) $}  \sigma_{2}^{-1}\sigma_{1}
\scalebox{0.7}{$-q^6 \left(q^8-q^6+q^4-q^2+1\right) $}  \sigma_{2}^{-1}\sigma_{1}^{-1}
 \displaybreak[1] \\  & \qquad 
\scalebox{0.7}{$+0$}  \sigma_{1}\sigma_{1}\sigma_{2}
\scalebox{0.7}{$+q^{16} \left(q^2-1\right) $}  \sigma_{1}\sigma_{1}\sigma_{2}^{-1}
 \displaybreak[1] \\  & \qquad 
\scalebox{0.7}{$+\left(3 q^{22}-6 q^{20}+4 q^{18}-q^{16}-2 q^{14}+5 q^{12}-3 q^{10}+3 q^8-q^4+q^2\right) $}  \sigma_{1}\sigma_{2}\sigma_{1}
 \displaybreak[1] \\  & \qquad 
\scalebox{0.7}{$+q^4 \left(q^{12}-2 q^{10}+q^8-q^6-q^4+q^2-1\right) $}  \sigma_{1}\sigma_{2}\sigma_{1}^{-1}
\scalebox{0.7}{$+q^{20} \left(q^2-1\right) $}  \sigma_{1}\sigma_{2}\sigma_{2}
\scalebox{0.7}{$-q^{12} $}  \sigma_{1}\sigma_{2}^{-1}\sigma_{1}
 \displaybreak[1] \\  & \qquad 
\scalebox{0.7}{$+0$}  \sigma_{1}\sigma_{2}^{-1}\sigma_{1}^{-1}
\scalebox{0.7}{$+q^4 \left(q^{12}-2 q^{10}+q^8-q^6-q^4+q^2-1\right) $}  \sigma_{1}^{-1}\sigma_{2}\sigma_{1}
 \displaybreak[1] \\  & \qquad 
\scalebox{0.7}{$+\left(q^{10}-q^8+q^6\right) $}  \sigma_{1}^{-1}\sigma_{2}\sigma_{1}^{-1}
\scalebox{0.7}{$+0$}  \sigma_{1}^{-1}\sigma_{2}\sigma_{2}
\scalebox{0.7}{$+0$}  \sigma_{1}^{-1}\sigma_{2}^{-1}\sigma_{1}
 \displaybreak[1] \\  & \qquad 
\scalebox{0.7}{$+q^8 \left(q^4-q^2+1\right) $}  \sigma_{1}^{-1}\sigma_{2}^{-1}\sigma_{1}^{-1}
\scalebox{0.7}{$+0$}  \sigma_{2}\sigma_{1}\sigma_{1}
\scalebox{0.7}{$+q^{20} \left(q^2-1\right) $}  \sigma_{2}\sigma_{2}\sigma_{1}
 \displaybreak[1] \\  & \qquad 
\scalebox{0.7}{$+0$}  \sigma_{2}\sigma_{2}\sigma_{1}^{-1}
\scalebox{0.7}{$+q^{16} \left(q^2-1\right) $}  \sigma_{2}^{-1}\sigma_{1}\sigma_{1}
\scalebox{0.7}{$+\left(q^{18}-q^{20}\right) $}  \sigma_{1}\sigma_{1}\sigma_{2}\sigma_{1}
 \displaybreak[1] \\  & \qquad 
\scalebox{0.7}{$+\left(q^{18}-q^{20}\right) $}  \sigma_{1}\sigma_{2}\sigma_{1}\sigma_{1}
\scalebox{0.7}{$+\left(q^{18}-q^{20}\right) $}  \sigma_{1}\sigma_{2}\sigma_{2}\sigma_{1}
\scalebox{0.7}{$+0$}  \sigma_{2}\sigma_{1}\sigma_{1}\sigma_{2}
 \\ & \Bigg) \\
\mathfig{0.1}{LY} =&  \left(\left(q^2-1\right)^2 \left(q^2+1\right) \left(q^4+1\right) \left(q^4-q^2+1\right)^2\right)^{-1} \times \Bigg( \\ & \qquad
\scalebox{0.7}{$-q^2 \left(q^2-1\right) \left(q^4+1\right) \left(q^{24}-q^{22}+2 q^{18}-2 q^{16}+q^{12}-2 q^{10}-q^4+q^2-1\right) $}  \mathbf{1}
 \displaybreak[1] \\  & \qquad 
\scalebox{0.7}{$-q^2 \left(q^{30}-3 q^{28}+5 q^{26}-4 q^{24}+2 q^{20}-3 q^{18}+q^{16}+2 q^{14}-q^{12}+2 q^{10}-2 q^6+2 q^4-2 q^2+1\right) $}  \sigma_{1}
 \displaybreak[1] \\  & \qquad 
\scalebox{0.7}{$+\left(-q^{26}+2 q^{24}-2 q^{22}+q^{20}+q^{18}-q^{16}+q^{14}-q^{10}+q^8-q^6+q^4\right) $}  \sigma_{1}^{-1}
 \displaybreak[1] \\  & \qquad 
\scalebox{0.7}{$+q^2 \left(q^{26}-2 q^{24}+2 q^{22}-4 q^{18}+8 q^{16}-9 q^{14}+5 q^{12}-5 q^8+6 q^6-5 q^4+3 q^2-1\right) $}  \sigma_{2}
 \displaybreak[1] \\  & \qquad 
\scalebox{0.7}{$-q^4 \left(q^{18}-q^{16}+q^{12}-2 q^{10}+q^6-q^4+q^2-1\right) $}  \sigma_{2}^{-1}
 \displaybreak[1] \\  & \qquad 
\scalebox{0.7}{$+\left(q^{30}-2 q^{28}+3 q^{26}-2 q^{24}-q^{18}+q^{12}\right) $}  \sigma_{1}\sigma_{1}
 \displaybreak[1] \\  & \qquad 
\scalebox{0.7}{$-q^2 \left(q^{24}-3 q^{22}+3 q^{20}+q^{18}-5 q^{16}+9 q^{14}-7 q^{12}+3 q^{10}-3 q^6+3 q^4-2 q^2+1\right) $}  \sigma_{1}\sigma_{2}
 \displaybreak[1] \\  & \qquad 
\scalebox{0.7}{$+\left(q^{18}-q^{10}+q^8-q^6+q^4\right) $}  \sigma_{1}\sigma_{2}^{-1}
 \displaybreak[1] \\  & \qquad 
\scalebox{0.7}{$+q^4 \left(q^4-q^2+1\right) \left(q^{14}-q^{12}+q^{10}-q^6+q^4-q^2+1\right) $}  \sigma_{1}^{-1}\sigma_{2}
\scalebox{0.7}{$-q^6 $}  \sigma_{1}^{-1}\sigma_{2}^{-1}
 \displaybreak[1] \\  & \qquad 
\scalebox{0.7}{$-q^6 \left(q^{20}-2 q^{18}+2 q^{16}-q^{12}+3 q^{10}-2 q^8+q^4-2 q^2+1\right) $}  \sigma_{2}\sigma_{1}
 \displaybreak[1] \\  & \qquad 
\scalebox{0.7}{$+q^6 \left(q^4-q^2+1\right) \left(q^8+q^2-1\right) $}  \sigma_{2}\sigma_{1}^{-1}
\scalebox{0.7}{$-q^{14} \left(q^2-1\right) \left(q^4+1\right) \left(q^6-q^2+1\right) $}  \sigma_{2}\sigma_{2}
 \displaybreak[1] \\  & \qquad 
\scalebox{0.7}{$-q^{12} \left(q^{10}-2 q^8+q^6-q^4-q^2+1\right) $}  \sigma_{2}^{-1}\sigma_{1}
 \displaybreak[1] \\  & \qquad 
\scalebox{0.7}{$-q^2 \left(q^{14}-q^{10}+q^8+q^6-2 q^4+2 q^2-1\right) $}  \sigma_{2}^{-1}\sigma_{1}^{-1}
\scalebox{0.7}{$+q^{14} \left(q^4-q^2+1\right) $}  \sigma_{1}\sigma_{1}\sigma_{2}
 \displaybreak[1] \\  & \qquad 
\scalebox{0.7}{$-q^{16} $}  \sigma_{1}\sigma_{1}\sigma_{2}^{-1}
\scalebox{0.7}{$+q^6 \left(2 q^{18}-4 q^{16}+2 q^{14}+q^{12}-3 q^{10}+5 q^8-2 q^6+q^4-1\right) $}  \sigma_{1}\sigma_{2}\sigma_{1}
 \displaybreak[1] \\  & \qquad 
\scalebox{0.7}{$+q^6 \left(q^4+1\right) \left(q^8-2 q^6+q^2-1\right) $}  \sigma_{1}\sigma_{2}\sigma_{1}^{-1}
\scalebox{0.7}{$+q^{14} \left(q^{10}-q^8+q^4-q^2+1\right) $}  \sigma_{1}\sigma_{2}\sigma_{2}
 \displaybreak[1] \\  & \qquad 
\scalebox{0.7}{$-q^{14} $}  \sigma_{1}\sigma_{2}^{-1}\sigma_{1}
\scalebox{0.7}{$+\left(q^2+1\right) \left(q^5-q^3+q\right)^2 $}  \sigma_{1}\sigma_{2}^{-1}\sigma_{1}^{-1}
\scalebox{0.7}{$-q^8 \left(q^2-1\right) \left(q^4+1\right) $}  \sigma_{1}^{-1}\sigma_{2}\sigma_{1}
 \displaybreak[1] \\  & \qquad 
\scalebox{0.7}{$+q^8 \left(q^4-q^2+1\right) $}  \sigma_{1}^{-1}\sigma_{2}\sigma_{1}^{-1}
\scalebox{0.7}{$-q^{16} \left(q^4-q^2+1\right) $}  \sigma_{1}^{-1}\sigma_{2}\sigma_{2}
\scalebox{0.7}{$+0$}  \sigma_{1}^{-1}\sigma_{2}^{-1}\sigma_{1}
 \displaybreak[1] \\  & \qquad 
\scalebox{0.7}{$+\left(-q^8+q^6-q^4\right) $}  \sigma_{1}^{-1}\sigma_{2}^{-1}\sigma_{1}^{-1}
\scalebox{0.7}{$+0$}  \sigma_{2}\sigma_{1}\sigma_{1}
\scalebox{0.7}{$+q^{22} \left(q^2-1\right) $}  \sigma_{2}\sigma_{2}\sigma_{1}
 \displaybreak[1] \\  & \qquad 
\scalebox{0.7}{$+0$}  \sigma_{2}\sigma_{2}\sigma_{1}^{-1}
\scalebox{0.7}{$+q^{18} \left(q^2-1\right) $}  \sigma_{2}^{-1}\sigma_{1}\sigma_{1}
\scalebox{0.7}{$+q^{12} \left(q^2-1\right) \left(q^4+1\right) $}  \sigma_{1}\sigma_{1}\sigma_{2}\sigma_{1}
 \displaybreak[1] \\  & \qquad 
\scalebox{0.7}{$+\left(q^{20}-q^{22}\right) $}  \sigma_{1}\sigma_{2}\sigma_{1}\sigma_{1}
\scalebox{0.7}{$+\left(q^{20}-q^{22}\right) $}  \sigma_{1}\sigma_{2}\sigma_{2}\sigma_{1}
\scalebox{0.7}{$+0$}  \sigma_{2}\sigma_{1}\sigma_{1}\sigma_{2}
\\ & \Bigg)
\end{align*}

\bibliographystyle{gtart}
\bibliography{bibliography/bibliography}

This paper is available online at \arxiv{0907.0256}, and at
\url{http://tqft.net/G2-surjectivity}.

\end{document}

%% file: article/article_preamble.tex
\input{\pathtotrunk preamble.tex}

\usepackage{breakurl}

\ifpdf
\usepackage[pdftex]{graphicx}
\else
\usepackage[dvips]{graphicx}
\fi

\usepackage{color}

%% file: preamble.tex

\usepackage{amsmath,amssymb,amsfonts}
\usepackage{ifpdf}
\usepackage{microtype}
\usepackage{pstricks,pst-node}

\ifpdf
\usepackage[pdftex,all,color]{xy}
\else
\usepackage[all,color]{xy}
\fi

\SelectTips{cm}{}

\vfuzz2pt 
\hfuzz2pt 

\ifpdf
\newcommand{\pathtodiagrams}{\pathtotrunk diagrams/pdf/}
\else
\newcommand{\pathtodiagrams}{\pathtotrunk diagrams/eps/}
\fi

\newcommand{\mathfig}[2]{{\hspace{-3pt}\begin{array}{c}%
  \raisebox{-2.5pt}{\includegraphics[width=#1\textwidth]{\pathtodiagrams #2}}%
\end{array}\hspace{-3pt}}}

\ifpdf
\usepackage[pdftex,plainpages=false,hypertexnames=false,pdfpagelabels]{hyperref}
\else
\usepackage[dvips,plainpages=false,hypertexnames=false]{hyperref}
\fi
\newcommand{\arxiv}[1]{\href{http://arxiv.org/abs/#1}{\tt arXiv:\nolinkurl{#1}}}
\newcommand{\doi}[1]{\href{http://dx.doi.org/#1}{{\tt DOI:#1}}}

\newcommand{\googlebooks}[1]{(preview at \href{http://books.google.com/books?id=#1}{google books})}

\def\RCS$#1: #2 ${\expandafter\def\csname RCS#1\endcsname{#2}}
\RCS$Revision: 1118 $ \RCS$Date: 2008-09-11 23:13:19 -0700 (Thu, 11 Sep
2008) $

\theoremstyle{plain}

\newtheorem{prop}{Proposition}[section]

\newtheorem{thm}[prop]{Theorem}
\newtheorem{lem}[prop]{Lemma}
\newtheorem{cor}[prop]{Corollary}
\newtheorem*{cor*}{Corollary}

\newtheorem*{defn*}{Definition}             
\newtheorem{question}{Question}
\newenvironment{rem}{\noindent\textsl{Remark.}}{}  
\numberwithin{equation}{section}

\newcounter{comment}
\newcommand{\noop}[1]{}

\def\clap#1{\hbox to 0pt{\hss#1\hss}}


\newcommand{\Complex}{\mathbb C}

\newcommand{\cA}{\mathcal{A}}
\newcommand{\cB}{\mathcal{B}}

\newcommand{\cR}{\mathcal{R}}

\newcommand{\cT}{\mathcal{T}}

\renewcommand{\imath}{\mathfrak{i}}
\renewcommand{\jmath}{\mathfrak{j}}

\newcommand{\tensor}{\otimes}

\newcommand{\End}[2]{\operatorname{End}_{#1}\left(#2\right)}

\newcommand{\Rep}{\operatorname{\textbf{Rep}}}


%% file: article/top_matter.tex
\title{The braid group surjects onto $G_2$ tensor space}

\author{Scott~Morrison}
\address{Microsoft Station Q, CNSI Building, UC Santa Barbara 93106-6105
}%
\email{scott@math.berkeley.edu} \urladdr{http://tqft.net/}

\primaryclass{
17B37 
}
\secondaryclass{
18D10 
17B10 
20G42 
}
\keywords{Braiding, G2, Spider
}